\documentclass{birkjour}
\usepackage[T1]{fontenc}
\usepackage{lmodern}
\usepackage{amsmath,amssymb,mathtools,microtype}
\usepackage{xcolor}
\usepackage[colorlinks=true,linkcolor=blue!45!black,citecolor=blue!45!black,urlcolor=blue!45!black]{hyperref}
\hypersetup{pdftitle={Alternative Clifford algebras: central decomposition and octonion quotients},pdfauthor={Adam Chapman},pdfsubject={Structure of alternative Clifford algebras and their centers}}
\copyrightinfo{}{}
\renewcommand{\subjclassname}{Mathematics Subject Classification (2020)}
\newtheorem{theorem}{Theorem}[section]
\newtheorem{proposition}[theorem]{Proposition}
\newtheorem{lemma}[theorem]{Lemma}
\newtheorem{corollary}[theorem]{Corollary}
\theoremstyle{definition}

\theoremstyle{remark}

\newcommand{\Calt}{C^{\mathrm{alt}}}
\newcommand{\Coct}{C^{\mathrm{oct}}}
\newcommand{\Mat}{\mathrm M}
\newcommand{\OO}{\mathcal O}
\newcommand{\Cent}{\mathrm Z}
\newcommand{\tr}{\operatorname{tr}}
\newcommand{\Tr}{\operatorname{Tr}}
\newcommand{\Spec}{\operatorname{Spec}}
\newcommand{\Frac}{\operatorname{Frac}}
\newcommand{\pair}[2]{\langle #1,#2\rangle}
\newcommand{\doi}[1]{\url{https://doi.org/#1}}
\newcommand{\revision}[1]{{\color{black}\hypersetup{allcolors=blue}#1}}
\numberwithin{equation}{section}
\title[Alternative Clifford algebras]{Alternative Clifford algebras:\\central decomposition and octonion quotients}
\author{Adam Chapman}
\address{School of Computer Science\\The Academic College of Tel Aviv-Yaffo\\Tel Aviv-Yafo, Israel}
\email{adamch@mta.ac.il}
\urladdr{https://orcid.org/0000-0001-6968-8139}
\date{28 October 2026}
\begin{document}
\begin{abstract}
We determine the structure of the alternative Clifford algebra of a nondegenerate quadratic form of dimension at least four over a field of characteristic different from two. A canonical central idempotent decomposes it into the ordinary Clifford algebra and an octonion algebra over a finitely generated commutative center. The octonion factor vanishes in dimension at least eight. In dimensions four through seven, its center is smooth of respective dimensions four, six, seven, and seven, and the full alternative Clifford algebra is semiprime. We give a universal Pl\"ucker presentation, an explicit affine quadric in dimension four, and a geometric description using orthogonal frames and the automorphism group of an octonion algebra. In characteristic two, the decomposition holds in polar rank at least four and the octonion factor vanishes in polar rank at least eight. Nonsingular forms of dimensions four and six have smooth centers, whereas regular seven-dimensional forms can have central nilpotents. As an application, Musgrave's octonion Clifford algebra is obtained by restricting the octonion factor to the reduced Zariski closure of the rational points of its center. It is infinite-dimensional over an algebraically closed field of characteristic different from two in dimensions four through seven.
\end{abstract}
\maketitle

\section{Introduction and main results}\label{sec:introduction}
The Clifford relations admit universal realizations in both associative and alternative algebras. The alternative Clifford algebra retains the same quadratic relations while allowing octonionic representations. Its structure therefore brings together ordinary Clifford algebras, composition algebras, and the geometry of quadratic forms. The problem addressed here is to describe this universal alternative algebra and its center in dimension at least four.

Musgrave studied representations of alternative Clifford algebras and introduced the octonion Clifford algebra by intersecting the kernels of their octonion-valued representations \cite{Musgrave}. Chapman and Vishne determined the alternative Clifford algebra of a nondegenerate ternary form, obtaining an octonion ring with polynomial center in one variable \cite{CV}. In higher dimension, a description of the full universal algebra requires identifying its associative and octonionic parts and determining whether any additional ideal remains. We prove that these two parts are complementary ideals defined by a canonical central idempotent.

The description has three parts. First, we establish the central decomposition and the dimension and polar-rank cutoffs for its octonion factor. Second, we present the center of that factor and describe its geometry through spaces of frames; in dimension four the presentation is a single affine quadric. Third, we extend the decomposition to characteristic two and apply it to Musgrave's octonion Clifford algebra, whose structure depends on rational points of the center.

The main ingredient is the coordinatization theorem of L\'opez Sol\'is and Shestakov for alternative algebras containing $\Mat_2(F)$ \cite{LS}, together with the functorial formulation of Grishkov and Shestakov \cite{GS}. In characteristic different from two, the commutator of two even generators satisfies a cubic identity. Its square yields a central idempotent, and the Clifford relations force all the idempotents obtained from different pairs to agree. In characteristic two, a suitable commutator is itself idempotent. Descent then gives the decomposition over the original field.

We work with unital algebras and unital homomorphisms. The field has characteristic different from two in Sections~\ref{sec:coordinates}--\ref{sec:geometry}; Section~\ref{sec:char-two} treats characteristic two, and Section~\ref{sec:octonion-clifford} states its characteristic assumptions explicitly. The comparisons in Section~\ref{sec:comparison} concern characteristic different from two. For a quadratic space $(V,q)$, write $A(q)=\Calt(q)$ for its alternative Clifford algebra and $C(q)$ for its ordinary Clifford algebra. \revision{We use $\mathcal A(D,M)$ for the coordinate algebra in \eqref{eq:coordinates} and $\mathcal R$ for the presentation of $R_q$ after normalization in Section~\ref{sec:presentation}. The variables $T_i$ in \eqref{eq:quadric} map to the central elements in \eqref{eq:T-def}.} In every characteristic the defining relation is $v^2=q(v)1$ for $v\in V$. If $\operatorname{char}F\ne2$ and $q=\langle a_1,\ldots,a_n\rangle$, the defining generators satisfy
\begin{equation}\label{eq:clifford}
 x_i^2=a_i,\qquad x_ix_j+x_jx_i=0\quad(i\ne j).
\end{equation}
There is a canonical surjection $\pi:A(q)\to C(q)$, whose kernel is the ideal generated by all associators.

An octonion algebra over a commutative ring means a unital composition algebra with nonsingular norm, locally free of rank eight. Its underlying module is not assumed to be free. The center of an alternative algebra consists of the elements that commute and associate with every pair of elements.

\begin{theorem}\label{thm:main}
Let $\operatorname{char}F\ne2$, and let $q$ be a nondegenerate quadratic form of dimension $n\ge4$ over $F$. There is a central idempotent $e\in A(q)$ such that $\pi$ restricts to an isomorphism $eA(q)\simeq C(q)$ and
\begin{equation}\label{eq:decomposition}
 A(q)\simeq C(q)\times\OO_q,\qquad \OO_q=(1-e)A(q).
\end{equation}
If $4\le n\le7$, then $\OO_q$ is a nonzero octonion algebra over its center $R_q$. The ring $R_q$ is a smooth finitely generated commutative $F$-algebra of pure dimension
\begin{equation}\label{eq:dimension}
 \dim R_q=7n-\frac{n(n+1)}2-14.
\end{equation}
It is geometrically integral for $4\le n\le6$ and has two geometric connected components for $n=7$.

If $n\ge8$, then $\OO_q=0$, so $\pi:A(q)\to C(q)$ is an isomorphism. In every dimension $n\ge4$, the algebra $A(q)$ is semiprime and its center is reduced.
\end{theorem}

For $q=\langle a_1,\ldots,a_n\rangle$, put
\[
 w_{ijk}=x_i(x_jx_k)+(x_jx_k)x_i.
\]
The idempotent in Theorem~\ref{thm:main} is
\begin{equation}\label{eq:intrinsic-e}
 e=-\frac{[w_{123},w_{124}]^2}{64a_1^2a_2^2a_3a_4}.
\end{equation}
The resulting decomposition is canonical, even though this expression uses an orthogonal basis.

In dimension four the center has a particularly short presentation.

\begin{theorem}\label{thm:four}
If $\operatorname{char}F\ne2$ and $q=\langle a,b,c,d\rangle$ with $abcd\ne0$, then
\begin{equation}\label{eq:quadric}
 R_q\simeq
 F[T_1,T_2,T_3,T_4,Z]\big/
 \left(Z^2-aT_1^2-bT_2^2-cT_3^2-dT_4^2-abcd\right).
\end{equation}
In particular, $\Cent(A(q))\simeq F\times R_q$.
\end{theorem}

The dimensions in \eqref{eq:dimension} and the vanishing of the octonionic quotient above dimension seven were announced in \cite[Theorem~8]{Vishne}. The structural contribution here is the direct-product decomposition of the full universal algebra, together with presentations, descent, and geometric proofs for the octonion factor and its center. In particular, the decomposition establishes semiprimeness. Section~\ref{sec:comparison} compares this consequence with the other assertions in that announcement under their precise characteristic and nondegeneracy hypotheses.

The characteristic-two statement is most naturally phrased using the polar form
\[
 b_q(v,w)=q(v+w)-q(v)-q(w).
\]
A quadratic form is called nonsingular when $b_q$ is nondegenerate; in characteristic two its dimension is then even.

\begin{theorem}\label{thm:char-two-main}
Suppose $\operatorname{char}F=2$ and $\operatorname{rank}b_q\ge4$. There is a canonical decomposition
\[
 A(q)\simeq C(q)\times\OO_q,
\]
where $\OO_q$ is either zero or an octonion algebra over a finitely generated commutative center $R_q$. If $\operatorname{rank}b_q\ge8$, then $\OO_q=0$.

If $q$ is nonsingular of dimension $n=4$ or $6$, then $\OO_q\ne0$ and $R_q$ is smooth and geometrically integral, of dimension $4$ or $7$, respectively. For every nonsingular $q$ of dimension at least four, $A(q)$ is semiprime and its center is reduced; for $n\ge8$, it equals $C(q)$.
\end{theorem}

In contrast, if $q=q_0\perp\langle a\rangle$, where $q_0$ is nonsingular of dimension six and $a\ne0$, then
\[
 R_q\simeq R_{q_0}\otimes_F F[t]/(t^2+a).
\]
When $a$ is a square, this center has a nonzero square-zero element and $A(q)$ is not semiprime. Proposition~\ref{prop:char-two-seven} proves this assertion and explains why the seven-dimensional conclusions of Theorem~\ref{thm:main} do not extend unchanged. \revision{For a regular five-dimensional form $q=q_0\perp\langle a\rangle$ with $q_0$ nonsingular of dimension four and $a\ne0$, the decomposition of Theorem~\ref{thm:char-two-main} still applies. We leave the description of $R_q$ in this case open.}

The decomposition also determines Musgrave's octonion Clifford algebra, denoted here by $\Coct_F(q)$. Its definition uses octonion algebras over $F$ itself, so rational points of the center matter. If $\OO_q\ne0$, put
\[
 I_F(R_q)=\bigcap_{\chi\in\operatorname{Hom}_{F\text{-alg}}(R_q,F)}\ker\chi,
\]
with the empty intersection interpreted as $R_q$.

\begin{corollary}\label{cor:octonion-main}
Under the hypotheses of Theorem~\ref{thm:main} or~\ref{thm:char-two-main},
\[
 \Coct_F(q)\simeq\OO_q\otimes_{R_q}R_q/I_F(R_q)
\]
when $\OO_q\ne0$, and $\Coct_F(q)=0$ when $\OO_q=0$. In characteristic two, the same intersection-of-kernels definition is used. If $F$ is algebraically closed of characteristic different from two and $4\le\dim q\le7$, then $\Coct_F(q)=\OO_q$ and is infinite-dimensional over $F$.
\end{corollary}

Section~\ref{sec:octonion-clifford} proves the corollary and distinguishes representations over $F$ from those over extension fields. The finite-dimensionality question in \cite[Question~3.4]{Musgrave} is discussed in Section~\ref{sec:comparison}.

Section~\ref{sec:coordinates} develops the coordinate tools. Section~\ref{sec:decomposition} proves the decomposition, its descent, and the dimension-eight cutoff in characteristic different from two. Section~\ref{sec:presentation} constructs the universal center and octonion factor. The explicit quadric is obtained in Section~\ref{sec:four}; the geometry of the centers and semiprimeness are established in Section~\ref{sec:geometry}. Section~\ref{sec:char-two} proves Theorem~\ref{thm:char-two-main}, gives an Artin--Schreier quadric presentation in dimension four, and treats the seven-dimensional exception. Section~\ref{sec:octonion-clifford} applies these results to the octonion Clifford algebra. Section~\ref{sec:comparison} concludes with the implications for earlier statements.

\section{Coordinatization and the Clifford relations}\label{sec:coordinates}
We use the following form of the coordinatization theorem \cite[Theorem~5.1]{LS}. An alternative algebra containing $\Mat_2(F)$ with the same identity has the form
\begin{equation}\label{eq:coordinates}
 \mathcal A(D,M)=\Mat_2(D)\oplus(M\oplus M),
\end{equation}
where $D$ is associative, $M$ is a left $D$-module with $[D,D]M=0$, and there is an alternating $D$-bilinear pairing
\[
 \pair{\ }{\ }:M\times M\longrightarrow\Cent(D)
\]
satisfying
\begin{equation}\label{eq:module-plucker}
 \pair{u}{v}w+\pair{v}{w}u+\pair{w}{u}v=0.
\end{equation}
The right action of $D$ on $M$ agrees with its left action. For $X,Y\in\Mat_2(D)$ and $u,v,w,z\in M$, multiplication is
\begin{equation}\label{eq:product}
 \begin{split}
 (X+(u,v))(Y+(w,z))={}&XY+
 \begin{pmatrix}
 -\pair{u}{z}&-\pair{v}{z}\\
 \pair{u}{w}&\pair{v}{w}
 \end{pmatrix}\\
 &+(w,z)X+(u,v)Y^*,
 \end{split}
\end{equation}
where
\[
 \begin{pmatrix}\alpha&\beta\\\gamma&\delta\end{pmatrix}^{*}
 =\begin{pmatrix}\delta&-\beta\\-\gamma&\alpha\end{pmatrix}.
\]
Conversely, these data define an alternative algebra. Homomorphisms of coordinate data induce homomorphisms fixing $\Mat_2(F)$; this is part of the categorical description in \cite{GS}.

Two elementary consequences will be used repeatedly.

\begin{lemma}\label{lem:center}
Let $J\subseteq D$ be the ideal generated by the pairing values. Then
\[
 [D,D]J=0.
\]
Every $h\in\Cent(D)$ gives a central element $hI_2\in\mathcal A(D,M)$.
\end{lemma}
\begin{proof}
For $d,d'\in D$ and $u,v\in M$, bilinearity gives
\[
 [d,d']\pair{u}{v}=\pair{[d,d']u}{v}=0.
\]
As the pairing values are central, the same holds for their multiples.

For the second assertion, a scalar matrix $hI_2$ commutes with $\Mat_2(D)$ and acts on either side of $(u,v)$ as $(hu,hv)$. It therefore commutes with the whole algebra. To verify associativity with it, it suffices by trilinearity to consider the matrix and module summands of \eqref{eq:coordinates}. With two matrix arguments, the assertion follows from associativity of $D$. With one module argument, it follows from the module law and centrality of $h$. With two module arguments, it follows from
\[
 \pair{hu}{v}=h\pair{u}{v}=\pair{u}{hv}
\]
in the matrix term of \eqref{eq:product}. Thus $(hI_2,x,y)=0$ for all $x,y$. Alternation of the associator gives the other two positions.
\end{proof}

\begin{lemma}\label{lem:local-octonion}
Suppose $D=R$ is commutative, $M$ is generated by finitely many elements $z_\alpha$, and their pairings generate the unit ideal of $R$. Then $M$ is locally free of rank two with perfect alternating pairing, and $\mathcal A(R,M)$ is an octonion algebra with center $R$.
\end{lemma}
\begin{proof}
On a principal open where $p=\pair{z_\alpha}{z_\beta}$ is invertible, \eqref{eq:module-plucker} gives
\[
 z=p^{-1}\bigl(\pair{z_\alpha}{z}z_\beta-
                 \pair{z_\beta}{z}z_\alpha\bigr)
 \qquad(z\in M).
\]
Thus $z_\alpha,z_\beta$ span. Pairing a linear relation between them with each of these two elements proves independence. After rescaling one basis vector, the pairing matrix is the standard symplectic matrix. Formula~\eqref{eq:product} then gives the split octonion algebra, as in \cite[Proposition~6.1]{LS}. The local conjugations and norms agree because they are given globally by
\begin{equation}\label{eq:norm}
 \overline{X+(u,v)}=X^*-(u,v),\qquad
 N(X+(u,v))=\det X+\pair{u}{v}.
\end{equation}
The center is $R$, as can be checked on the same open cover.
\end{proof}

We first normalize the first two coefficients of the quadratic form:
\[
 q=\langle1,1,b_1,\ldots,b_m\rangle,\qquad
 m=n-2,\qquad b_r\ne0.
\]
The first two generators give a copy of $\Mat_2(F)$. Indeed, they generate a homomorphic image of the binary Clifford algebra $\Mat_2(F)$, and its composition with the ordinary Clifford quotient is injective. Choose
\begin{equation}\label{eq:matrices}
 i=\begin{pmatrix}1&0\\0&-1\end{pmatrix},\qquad
 j=\begin{pmatrix}0&1\\1&0\end{pmatrix},\qquad
 k=ij=\begin{pmatrix}0&1\\-1&0\end{pmatrix}.
\end{equation}
\revision{For $X=\left(\begin{smallmatrix}a&b\\ c&d\end{smallmatrix}\right)\in\Mat_2(D)$,
\[
 iX+Xi=\begin{pmatrix}2a&0\\0&-2d\end{pmatrix},\qquad
 jX+Xj=\begin{pmatrix}b+c&a+d\\a+d&b+c\end{pmatrix}.
\]}
An element anticommuting with $i$ and $j$ has matrix component $sk$, and its module component is unrestricted. Hence in the coordinates of $A(q)$,
\begin{equation}\label{eq:generators}
 x_{r+2}=s_rk+(u_r,v_r)\qquad(1\le r\le m).
\end{equation}

\begin{lemma}\label{lem:universal-coordinates}
The remaining Clifford relations are equivalent to
\begin{equation}\label{eq:coordinate-relations}
 \begin{aligned}
 s_r^2+\pair{u_r}{v_r}&=-b_r,\\
 s_rs_t+s_ts_r+\pair{u_r}{v_t}+\pair{u_t}{v_r}&=0
 \qquad(r\ne t).
 \end{aligned}
\end{equation}
Moreover, $D$ is generated by the $s_r$ and the pairings of the $u_r,v_r$, and these $2m$ vectors generate $M$ as a $D$-module.
\end{lemma}
\begin{proof}
The square of $(u,v)$ is $-\pair{u}{v}I_2$. Its anticommutator with $sk$ is zero, while the anticommutator of $(u_r,v_r)$ and $(u_t,v_t)$ is
$-(\pair{u_r}{v_t}+\pair{u_t}{v_r})I_2$. This proves \eqref{eq:coordinate-relations}.

The coordinate generators can be recovered inside the algebra. For $x=sk+(u,v)$,
\begin{equation}\label{eq:recovery}
 sI_2=-\tfrac12(kx+xk),\qquad (u,v)=x-sk.
\end{equation}
Multiplication by the matrix units separates the two coordinates and transfers them between the two positions. Their products recover the pairing values through \eqref{eq:product}.

More formally, let $D_0$ be the subalgebra generated by the $s_r$ and those pairing values, and let $M_0$ be the $D_0$-span of the $u_r,v_r$. Bilinearity shows that $\mathcal A(D_0,M_0)$ is a subalgebra containing all the Clifford generators. It is therefore the whole of $A(q)$, giving $D_0=D$ and $M_0=M$. Conversely, any coordinate data satisfying \eqref{eq:coordinate-relations} give a representation of the Clifford generators by \eqref{eq:generators}. This observation will justify the universal constructions below.
\end{proof}

\section{The central decomposition}\label{sec:decomposition}
Continue in the normalized case, with $m\ge2$. Let $J\subseteq D$ be the pairing ideal and, for $r\ne t$, set
\[
 c_{rt}=[s_r,s_t].
\]
Relations~\eqref{eq:coordinate-relations} imply that $s_r^2$, $s_t^2$, and $s_rs_t+s_ts_r$ are central in $D$. Thus
\begin{equation}\label{eq:commutator-square}
 c_{rt}^2=(s_rs_t+s_ts_r)^2-4s_r^2s_t^2\in\Cent(D).
\end{equation}
Modulo $J$, this square is $-4b_rb_t$. Since $c_{rt}J=0$ by Lemma~\ref{lem:center}, we obtain
\begin{equation}\label{eq:cubic}
 c_{rt}^3=-4b_rb_t\,c_{rt}.
\end{equation}

\begin{proposition}\label{prop:idempotent}
The elements
\[
 e_{rt}=-\frac{c_{rt}^2}{4b_rb_t}\qquad(r\ne t)
\]
are equal central idempotents of $A(q)$. For their common value $e$,
\[
 eM=eJ=0,\qquad ec_{rt}=c_{rt}.
\]
The ordinary Clifford quotient induces an isomorphism $eA(q)\simeq C(q)$.
\end{proposition}
\begin{proof}
Equation~\eqref{eq:cubic} gives $c_{rt}^4=-4b_rb_t\,c_{rt}^2$, hence $e_{rt}^2=e_{rt}$ and $e_{rt}c_{rt}=c_{rt}$. Each $c_{rt}$ annihilates $M$ and $J$, so the same is true of $e_{rt}$. Centrality in the whole alternative algebra follows from \eqref{eq:commutator-square} and Lemma~\ref{lem:center}.

In the $e_{rt}$-component the pairings vanish. Consequently its coordinate generators satisfy
\[
 (e_{rt}s_\ell)^2=-b_\ell e_{rt},\qquad
 (e_{rt}s_\ell)(e_{rt}s_h)+(e_{rt}s_h)(e_{rt}s_\ell)=0
 \quad(\ell\ne h).
\]
For any distinct $\ell,h$, the square of their commutator is
$-4b_\ell b_h e_{rt}$. Hence $e_{rt}e_{\ell h}=e_{rt}$.
Interchanging the two pairs gives $e_{rt}e_{\ell h}=e_{\ell h}$, proving equality.

Since $eM=0$, the algebra $eA(q)=\Mat_2(eD)$ is associative. In the ordinary Clifford quotient the Cayley component is zero and the even coordinates satisfy the same Clifford relations, so $\pi(e)=1$. There are therefore homomorphisms
\[
 C(q)\longrightarrow eA(q)\longrightarrow C(q)
\]
whose composite fixes the Clifford generators. The first map is surjective because the elements $ex_i$ generate $eA(q)$ with identity $e$. The maps are inverse.
\end{proof}

\begin{proposition}\label{prop:complement}
Let $f=1-e$. The algebra $R=fD$ is commutative and finitely generated over $F$. Either $fA(q)=0$, or $fA(q)$ is an octonion algebra over $R$, with center $R$.
\end{proposition}
\begin{proof}
Since $ec_{rt}=c_{rt}$, all the generators $fs_r$ commute. The remaining generators of $fD$ are pairing values, which are central. Thus $R=fD$ is commutative and finitely generated by Lemma~\ref{lem:universal-coordinates}.

The pairing ideal $fJ$ is the unit ideal of $R$. In fact, in $R/fJ$ we have
\[
 s_1^2=-b_1,\qquad s_2^2=-b_2,\qquad 2s_1s_2=0.
\]
The first two elements are units, so the third equality forces the quotient ring to be zero. Lemma~\ref{lem:local-octonion} now applies to $fM$, which is generated by the images of the $u_r,v_r$.
\end{proof}

We next descend the decomposition. Let $q=\langle a_1,\ldots,a_n\rangle$ be arbitrary and put $E=F(\sqrt{a_1},\sqrt{a_2})$. This is a finite separable extension. Formation of $A(q)$ commutes with this extension, \revision{since the free alternative algebra on $V$ base changes to that on $V_E$, and the square relations on a basis and their pairwise polarizations generate the defining ideal over either field.} Write
\[
 y_1=x_1/\sqrt{a_1},\qquad y_2=x_2/\sqrt{a_2},
 \qquad y_r=x_r\quad(r\ge3).
\]
For these normalized generators, \eqref{eq:product} gives
\[
 y_1(y_2y_{r+2})+(y_2y_{r+2})y_1=-2s_rI_2.
\]
Therefore the normalized idempotent is precisely the scalar extension of \eqref{eq:intrinsic-e}. Since that expression is defined over $F$, centrality and idempotence descend. The restriction of $\pi$ to $eA(q)$ becomes an isomorphism over $E$, and hence is already an isomorphism over $F$.

Put $B=(1-e)A(q)$ and $R_q=\Cent(B)$. Centers commute with this finite field extension: write an element of $B\otimes_F E$ in an $F$-basis of $E$ and compare coefficients in its commutators and associators with elements of $B$. Thus
\[
 \Cent(B\otimes_F E)=R_q\otimes_F E.
\]
Proposition~\ref{prop:complement} shows that this is the coefficient ring of the normalized octonion factor. Finite generation of $R_q$ and finite local freeness of $B$ of rank eight descend under the faithfully flat map $R_q\to R_q\otimes_F E$ whenever $B\ne0$.

For completeness, the composition structure descends as well. The map $x_i\mapsto-x_i$ extends to an anti-involution of $A(q)$: the opposite of an alternative algebra is alternative, and these images satisfy the same defining relations. It fixes $e$ and, over $E$, restricts on $B$ to the conjugation in \eqref{eq:norm}, because the two anti-involutions agree on the generators. Consequently
\[
 b+\bar b\in R_q,\qquad b\bar b\in R_q
\]
by faithful descent. The map $N(b)=b\bar b$ is a quadratic, multiplicative, nonsingular norm, since these properties hold after extension to $E$. This proves the decomposition over $F$.

The decomposition is canonical. Indeed, $B=\ker\pi$, and $eA(q)$ is its two-sided annihilator: the central product decomposition shows that an element annihilating $B$ has zero $B$-component, by multiplication with the identity $f$ of $B$. Thus $e$ is the identity of an intrinsically defined ideal.

\begin{lemma}\label{lem:seven}
Let $O$ be an octonion algebra over a field $L$ of characteristic different from two. If $n\ge2$ and elements $y_1,\ldots,y_n\in O$ satisfy
\[
 y_i^2=a_i\in L^\times,\qquad
 y_iy_j+y_jy_i=0\quad(i\ne j),
\]
then $n\le7$.
\end{lemma}
\begin{proof}
Every $y_i$ is invertible and nonscalar, since a nonzero scalar cannot anticommute with another invertible element. The quadratic identity
\[
 y_i^2-\Tr(y_i)y_i+N(y_i)=0
\]
then implies $\Tr(y_i)=0$ and $N(y_i)=-a_i$. For pure octonions,
$b_N(y_i,y_j)=-(y_iy_j+y_jy_i)$. Hence their Gram matrix is
$\operatorname{diag}(-2a_1,\ldots,-2a_n)$, and they are linearly independent in the seven-dimensional trace-zero subspace.
\end{proof}

\begin{corollary}\label{cor:eight}
If $\dim q\ge8$, then $\pi:A(q)\to C(q)$ is an isomorphism. In particular, $\dim_F A(q)=2^{\dim q}$.
\end{corollary}
\begin{proof}
If $R_q\ne0$, choose a maximal ideal of $R_q$. The fiber of $B$ is an octonion algebra over its residue field, and the images of the Clifford generators contradict Lemma~\ref{lem:seven}. Thus $B=0$.
\end{proof}

The argument uses more than the absence of octonionic field-valued representations: Proposition~\ref{prop:complement} identifies the entire complementary ideal as an octonion algebra over its center. This is what excludes a residual nonassociative radical.

\section{The universal center and octonion factor}\label{sec:presentation}
We present the center for $q=\langle1,1,b_1,\ldots,b_m\rangle$, where $m\ge2$.
Introduce variables $s_1,\ldots,s_m$ and $p_{\alpha\beta}$ for
$1\le\alpha<\beta\le2m$, and put
$p_{\beta\alpha}=-p_{\alpha\beta}$, $p_{\alpha\alpha}=0$.
Let $\mathcal R$ be the commutative polynomial algebra on these variables modulo
\begin{equation}\label{eq:center-relations}
 \begin{aligned}
 p_{\alpha\beta}p_{\gamma\delta}
 -p_{\alpha\gamma}p_{\beta\delta}
 +p_{\alpha\delta}p_{\beta\gamma}&=0,\\
 p_{2r-1,2r}&=-(s_r^2+b_r),\\
 p_{2r-1,2t}+p_{2t-1,2r}&=-2s_rs_t\quad(r\ne t).
 \end{aligned}
\end{equation}
The first family is imposed for all indices. \revision{Over an algebraic closure $K$, the torsor \eqref{eq:sl2-torsor} identifies $R_q\otimes_F K$ with $K[X_m]^{\mathrm{SL}_2}$. For the vectors in $K^2$, the first fundamental theorem for $\mathrm{SL}_2$ gives the determinant generators $p_{\alpha\beta}$, and the second gives their Pl\"ucker relations in the first family.} Let $P$ be the $\mathcal R$-module on $z_1,\ldots,z_{2m}$ with relations
\begin{equation}\label{eq:P-relations}
 p_{\alpha\beta}z_\gamma+
 p_{\beta\gamma}z_\alpha+
 p_{\gamma\alpha}z_\beta=0.
\end{equation}

\begin{proposition}\label{prop:universal}
The pairing $\pair{z_\alpha}{z_\beta}=p_{\alpha\beta}$ is well defined on $P$. If $\mathcal R\ne0$, it is perfect and $P$ is locally free of rank two. Moreover,
\begin{equation}\label{eq:universal-model}
 R_q\simeq\mathcal R,\qquad
 \OO_q\simeq\mathcal A(\mathcal R,P).
\end{equation}
\end{proposition}
\begin{proof}
Pairing a defining relation in \eqref{eq:P-relations} with $z_\delta$ gives a Pl\"ucker relation in \eqref{eq:center-relations}. Thus the pairing is well defined and satisfies \eqref{eq:module-plucker}. Modulo all the $p_{\alpha\beta}$, the last two families of relations give $s_1^2=-b_1$, $s_2^2=-b_2$, and $2s_1s_2=0$. Therefore the $p_{\alpha\beta}$ generate the unit ideal. Lemma~\ref{lem:local-octonion} applies, so
\[
 B_0=\mathcal A(\mathcal R,P)
\]
is an octonion algebra with center $\mathcal R$, unless $\mathcal R=0$.

Take in $B_0$ the constant matrices $i,j$ and the elements
\[
 s_rk+(z_{2r-1},z_{2r})\qquad(1\le r\le m).
\]
By Lemma~\ref{lem:universal-coordinates}, they satisfy the Clifford relations. Universality gives $A(q)\to B_0$. The even coordinates commute, so this map kills $e$ and factors through $\OO_q$.

Conversely, in the coordinates of $\OO_q$, the even generators commute, their pairings satisfy the Pl\"ucker relations, and \eqref{eq:coordinate-relations} holds. They define a homomorphism $\mathcal R\to R_q$ and a compatible module map $P\to fM$. Formula~\eqref{eq:product} gives a homomorphism $B_0\to\OO_q$. The two composites fix the constant matrices, the even coordinates, and the module generators. These generate the respective algebras by \eqref{eq:recovery}. The maps are inverse, and their restrictions identify the centers. The same argument covers the zero-ring case.
\end{proof}

This construction establishes the full universal property of the octonion factor. In particular, the ring in \eqref{eq:center-relations} is not merely the coordinate ring of a family of representations. For an arbitrary diagonal form, the presentation holds after extending to $F(\sqrt{a_1},\sqrt{a_2})$, and the preceding section gives its descent.

\section{The affine quadric in dimension four}\label{sec:four}
Let first $q=\langle1,1,c,d\rangle$. Write the two even coordinates as $s,t$ and the four module generators as $z_1,\ldots,z_4$. In the commutative octonion factor, put
\[
 u=p_{13},\qquad v=p_{24},\qquad h=p_{14}+st.
\]
Then
\begin{equation}\label{eq:four-p}
 \begin{aligned}
 p_{12}&=-(s^2+c),&p_{34}&=-(t^2+d),\\
 p_{13}&=u,&p_{24}&=v,\\
 p_{14}&=h-st,&p_{23}&=h+st.
 \end{aligned}
\end{equation}
There is one independent Pl\"ucker relation,
\[
 p_{12}p_{34}-p_{13}p_{24}+p_{14}p_{23}=0.
\]
Substitution gives
\begin{equation}\label{eq:split-quadric}
 R_q\simeq F[s,t,h,u,v]/(uv-h^2-ds^2-ct^2-cd).
\end{equation}
The octonion algebra is $\mathcal A(R_q,P)$, where $P$ is given by
\eqref{eq:P-relations} and \eqref{eq:four-p}.

To express the center over the original field, let the $x_i$ denote their images in the octonion factor and define
\begin{equation}\label{eq:T-def}
 T_1=\tfrac12w_{234},\quad T_2=\tfrac12w_{134},\quad
 T_3=\tfrac12w_{124},\quad T_4=\tfrac12w_{123},
\end{equation}
and
\begin{equation}\label{eq:Z-def}
 Z=\tfrac12\bigl((x_1x_2)(x_3x_4)+(x_3x_4)(x_1x_2)\bigr).
\end{equation}
For the normalized form, \eqref{eq:product} gives
\begin{equation}\label{eq:change}
 T_1=\frac{u-v}{2},\qquad T_2=-h,\qquad
 T_3=-t,\qquad T_4=-s,\qquad Z=\frac{u+v}{2}.
\end{equation}
For example, the matrix component of $x_3x_4$ is
\[
 X=-stI_2+\begin{pmatrix}-p_{14}&-p_{24}\\p_{13}&p_{23}\end{pmatrix},
 \qquad \tr X=0.
\]
Its symmetrized products with $j,i,k$ are respectively
$\tr(jX)I_2$, $\tr(iX)I_2$, and $\tr(kX)I_2$; the module terms cancel. These traces give $2T_1,2T_2,2Z$. The remaining two identities follow from the calculation of $w_{123},w_{124}$ in the descent argument.

The inverse change of coordinates is
\[
 u=Z+T_1,\quad v=Z-T_1,\quad h=-T_2,\quad t=-T_3,\quad s=-T_4.
\]
Equation~\eqref{eq:split-quadric} becomes
\begin{equation}\label{eq:normalized-quadric}
 Z^2-T_1^2-T_2^2-cT_3^2-dT_4^2=cd.
\end{equation}

Now take $q=\langle a,b,c,d\rangle$ and work over $E=F(\sqrt a,\sqrt b)$. For the normalized generators, the elements in \eqref{eq:T-def}--\eqref{eq:Z-def} satisfy
\[
 T'_1=\frac{T_1}{\sqrt b},\quad
 T'_2=\frac{T_2}{\sqrt a},\quad
 T'_3=\frac{T_3}{\sqrt{ab}},\quad
 T'_4=\frac{T_4}{\sqrt{ab}},\quad
 Z'=\frac{Z}{\sqrt{ab}}.
\]
Multiplication of \eqref{eq:normalized-quadric} by $ab$ yields \eqref{eq:quadric}. The elements $T_i,Z$ are central over $F$, since they are central over $E$. The resulting homomorphism from the ring displayed in \eqref{eq:quadric} to $R_q$ becomes an isomorphism after extension to $E$ and hence is an isomorphism over $F$. Since $C(q)$ is central simple in dimension four, its center is $F$. This proves Theorem~\ref{thm:four}.

The quadric is smooth: the gradient can vanish only at the origin, which is not on the quadric because $abcd\ne0$. Its quadratic part has rank five over an algebraic closure, whereas a product of two linear forms has rank at most two. Thus its defining polynomial is geometrically irreducible. This also proves directly that $R_q$ is geometrically integral of dimension four.

\section{Orthogonal frames and the geometry of the centers}\label{sec:geometry}
We first prove the dimension and smoothness assertions using the explicit presentation in Section~\ref{sec:presentation}.

\begin{proposition}\label{prop:geometry}
For $4\le n\le7$, the ring $R_q$ is nonzero and smooth of pure dimension \eqref{eq:dimension}. It is geometrically integral for $4\le n\le6$; for $n=7$, it has exactly two geometric connected components, each geometrically integral of dimension seven.
\end{proposition}
\begin{proof}
Work over an algebraic closure $K$ of $F$ and normalize the first two coefficients. Let $U=K^2$ with its standard alternating form $\omega$, and define the nonsingular five-dimensional quadratic space
\[
 W=K\oplus U\oplus U,\qquad Q(s,u,v)=s^2+\omega(u,v).
\]
Let $X_m$ be the scheme of ordered orthogonal $m$-frames in $W$, with prescribed norms $-b_1,\ldots,-b_m$, where $m=n-2$.

There is a morphism
\[
 X_m\longrightarrow\Spec\mathcal R
\]
sending a frame $(s_r,u_r,v_r)_r$ to the coordinates $s_r$ and the pairings of the ordered list $u_1,v_1,\ldots,u_m,v_m$. The norm and orthogonality equations are exactly the last two families in \eqref{eq:center-relations}, and the pairings satisfy the Pl\"ucker relations.

This morphism is the symplectic frame bundle of $P$. Indeed, for any commutative $K$-algebra $S$, a symplectic trivialization
$P\otimes_{\mathcal R}S\simeq S^2$ sends the $z_\alpha$ to vectors whose pairings are the prescribed $p_{\alpha\beta}$, producing an $S$-point of $X_m$. Conversely, those vectors satisfy \eqref{eq:P-relations} and define a pairing-preserving map $P\otimes S\to S^2$. Wherever $p_{\alpha\beta}$ is invertible, the two corresponding vectors form a basis, so this map is an isomorphism. These opens cover the base by Proposition~\ref{prop:universal}. The constructions are inverse and commute with base change. Hence
\begin{equation}\label{eq:sl2-torsor}
 X_m\longrightarrow\Spec\mathcal R
 \quad\hbox{is a Zariski locally trivial $\mathrm{SL}_2$-torsor.}
\end{equation}

For $m\le5$, the frame scheme is nonempty, smooth, and of dimension
\[
 \dim X_m=5m-\frac{m(m+1)}2.
\]
To see smoothness directly, the Gram matrix of the frame is nonsingular. Varying its vectors within their span realizes every infinitesimal symmetric variation of the Gram matrix, since $2$ is invertible. Thus the $m(m+1)/2$ Gram equations have independent differentials. Equivalently, after choosing a reference frame,
$X_m\simeq\mathrm O_5/\mathrm O_{5-m}$.

For $m<5$, the group $\mathrm{SO}_5$ acts transitively: an isometry extending a frame can have its determinant adjusted on the nonzero orthogonal complement. Thus $X_m$ is geometrically integral. For $m=5$, it is $\mathrm O_5$, which has two geometrically integral components. These facts use $\operatorname{char}K\ne2$.

Smoothness descends along \eqref{eq:sl2-torsor} \revision{(Stacks Project, \href{https://stacks.math.columbia.edu/tag/02K5}{Tag~02K5})}, and subtracting $\dim\mathrm{SL}_2=3$ gives
\[
 \dim\mathcal R
 =5m-\frac{m(m+1)}2-3
 =7n-\frac{n(n+1)}2-14.
\]
For $m<5$, the surjective image of the integral scheme $X_m$ is irreducible; smoothness makes it reduced. For $m=5$, the connected group $\mathrm{SL}_2$ preserves each orientation component. The two components therefore descend to disjoint open and closed components of the quotient, each integral. This proves all the geometric assertions and nonvanishing.
\end{proof}

Thus the center dimensions are
\[
\begin{array}{c|cccc}
 n&4&5&6&7\\ \hline
 \dim R_q&4&6&7&7.
\end{array}
\]
For $n\ge8$ the octonion factor, and hence its center, is zero. Its absence should not be interpreted as a nonzero zero-dimensional parameter space.

We next give an intrinsic interpretation. For an octonion algebra $O$, let $O_0$ be its trace-zero subspace. We use the standard facts that the automorphism group of the split octonion algebra is the smooth group of type $G_2$, of dimension fourteen, and that octonion algebras become split in the faithfully flat topology; see \cite{SV} for the field theory and \cite{AG} for the group-scheme formulation over rings.

\begin{lemma}\label{lem:frame-generation}
Let $S$ be a commutative $F$-algebra, and let $O$ be an octonion algebra over $S$. If $4\le n\le7$ and $y_1,\ldots,y_n\in O_0$ satisfy \eqref{eq:clifford} with $a_i\in F^\times$, then they generate $O$ as an $S$-algebra.
\end{lemma}
\begin{proof}
Generation can be checked after a faithfully flat extension. Normalize the first two coefficients and split the quaternion subalgebra generated by $y_1,y_2$. Locally the octonion algebra has the form $\mathcal A(S,M)$ with $M$ a symplectic module of rank two. Relative to that quaternion algebra, write the remaining elements as in \eqref{eq:generators}, with $s_r\in S$. Their pairings generate the unit ideal, by the same argument with $s_1,s_2$ used in Proposition~\ref{prop:complement}. On a cover where one of these pairings is invertible, the two corresponding module vectors form a basis of $M$. The recovery operations \eqref{eq:recovery} and the matrix units show that the algebra generated by the $y_i$ contains both copies of $M$ as well as the quaternion algebra. It is therefore all of $O$ locally and hence globally.
\end{proof}

\begin{proposition}\label{prop:g2}
Over an algebraic closure $K$ of $F$, let $O$ be a fixed split octonion algebra and let $Y_n$ be the scheme of isometric embeddings
\[
 (V_K,-q_K)\hookrightarrow(O_0,N).
\]
For $4\le n\le7$, the morphism
\[
 Y_n\longrightarrow\Spec(R_q\otimes_F K)
\]
is a $G_2$-torsor. In particular, its base is the faithfully flat sheaf quotient by $G_2$. For $n=6,7$, respectively, the base is isomorphic to
\[
 \mathrm{SO}_7/G_2,\qquad \mathrm O_7/G_2.
\]
\end{proposition}
\begin{proof}
A frame over a $K$-algebra $S$ yields a homomorphism
$A(q)\otimes_F K\to O\otimes_K S$. The elements $w_{ijk}$ are scalar in an octonion algebra when the generators are pure: $x_jx_k$ is pure for distinct orthogonal generators, and its anticommutator with $x_i$ is scalar. Thus \eqref{eq:intrinsic-e} maps to zero and the homomorphism factors through $\OO_q$.

By Lemma~\ref{lem:frame-generation}, its $S$-linear extension onto $O\otimes_K S$ is surjective. The center of $\OO_q$ consequently maps into the center $S$, giving a map $R_q\otimes_F K\to S$. The induced homomorphism
\[
 \OO_q\otimes_{R_q}S\longrightarrow O\otimes_K S
\]
is a surjection between locally free modules of rank eight, hence an isomorphism. Conversely, a map of centers and such an isomorphism give the frame by taking the universal generators. The constructions are functorial and inverse. Therefore $Y_n$ is the torsor of octonion trivializations of the universal factor, and its structure group is $G_2$.

For $n=7$, the full frame scheme is an orthogonal-group torsor, identified over $K$ with $\mathrm O_7$. For $n=6$, a prescribed orientation gives a unique completion to a full orthogonal frame, identifying the frame scheme with $\mathrm{SO}_7$. The automorphism group $G_2$ lies in $\mathrm{SO}_7$, since it is connected. Taking the quotients and using inversion to pass from left cosets to right cosets proves the stated identifications. \revision{For $n=4,5$, the pointwise stabilizer of a reference $n$-frame in $\mathrm{SO}_7$ is $\mathrm{SO}_{7-n}$ on its orthogonal complement. Thus the base is the double quotient $\mathrm{SO}_{7-n}\backslash\mathrm{SO}_7/G_2$ as a faithfully flat sheaf. For $n=4$ it is the affine quadric \eqref{eq:quadric} over $K$. For $n=7$, fix a reference isometry $\iota_0:(V_K,-q_K)\to(O_0,N)$. The sign $\det(\iota_0^{-1}\iota)$ of a representation $\iota$ distinguishes the two components. It is independent of the octonion trivialization because $G_2\subset\mathrm{SO}_7$, and reversing the reference orientation interchanges the labels.}
\end{proof}

Over $F$, these parameter spaces can be twisted; their nonvanishing does not assert the existence of an $F$-rational specialization. The $G_2$ description gives a second dimension calculation:
the orthogonal $n$-frame scheme in dimension seven has dimension
$7n-n(n+1)/2$, and the torsor has fibers of dimension fourteen.

\begin{corollary}\label{cor:semiprime}
For every nondegenerate $q$ of dimension at least four, $A(q)$ is semiprime and $\Cent(A(q))$ is reduced.
\end{corollary}
\begin{proof}
The ordinary Clifford algebra is semisimple for a nondegenerate form in characteristic different from two, and its center is finite \'{e}tale over $F$. For $4\le n\le7$, Proposition~\ref{prop:geometry} makes $R_q$ a reduced noetherian ring, and $\OO_q$ is projective over $R_q$. Consequently it embeds into the product of its fibers over the fraction fields of the minimal-prime quotients of $R_q$. Each such fiber is a simple octonion algebra. A nilpotent ideal maps to zero in every one of these fibers and therefore is zero. Thus $\OO_q$ is semiprime; its center $R_q$ is reduced. The decomposition \eqref{eq:decomposition} proves the assertion, with the case $n\ge8$ already covered by Corollary~\ref{cor:eight}.
\end{proof}

Propositions~\ref{prop:idempotent}--\ref{prop:complement}, the descent argument, Corollary~\ref{cor:eight}, Proposition~\ref{prop:geometry}, and Corollary~\ref{cor:semiprime} complete the proof of Theorem~\ref{thm:main}.

\section{Characteristic two}\label{sec:char-two}
Throughout this section, $\operatorname{char}F=2$. Write $\beta=b_q$ for the polar form of $q$. The results of \cite{LS,GS}, and Lemmas~\ref{lem:center} and~\ref{lem:local-octonion}, hold in this characteristic as well. \revision{The pairing in \cite[Theorem~5.1]{LS} is stated to be skew-symmetric. Setting $w=u$ in \eqref{eq:module-plucker} gives $\pair{u}{u}M=0$. By bilinearity, $\pair{u}{u}$ annihilates every pairing value, so $\pair{u}{u}=0$ wherever the pairing ideal is the unit ideal, as in Lemma~\ref{lem:local-octonion} and on the octonion factor. Even before passing to this factor, $m=(u,0)$ satisfies $mE_{11}=0$ by \eqref{eq:product}, so $\pair{u}{u}E_{21}=m^2E_{11}=m(mE_{11})=0$.} The ternary characteristic-two case already occurs in \cite[Corollary~4.6]{CV}; the issue addressed below is the decomposition in higher polar rank.

\subsection{The idempotent and the rank cutoff}
Assume $\operatorname{rank}\beta\ge4$. Choose a nonsingular binary subform $q_1$ and write $q=q_1\perp q'$. After a finite separable extension, $q_1$ is hyperbolic. We first work over such a field and write $q=H\perp q'$. Choose generators $i,j$ of $H$ with
\[
 i^2=j^2=0,\qquad ij+ji=1.
\]
They generate $\Mat_2(F)$, with $i=E_{12}$ and $j=E_{21}$. As before, injectivity of this quaternion subalgebra follows by composition with the ordinary Clifford quotient.

Let $w_1,\ldots,w_m$ be a basis of the space of $q'$, put
\[
 a_r=q'(w_r),\qquad \beta_{rt}=\beta(w_r,w_t),\qquad m=\dim q-2,
\]
and choose the first two basis vectors so that $\beta_{12}=1$. The corresponding Clifford generators commute with $i$ and $j$. In the coordinatization \eqref{eq:coordinates}, this says precisely that they have the form
\begin{equation}\label{eq:char-two-generators}
 x_{r+2}=s_rI_2+(u_r,v_r).
\end{equation}
Indeed, their matrix parts centralize $\Mat_2(F)$, while a Cayley-module element commutes with every trace-zero matrix in characteristic two. Formula~\eqref{eq:product} now gives the relations
\begin{equation}\label{eq:char-two-relations}
 \begin{aligned}
 s_r^2+\pair{u_r}{v_r}&=a_r,\\
 [s_r,s_t]+\pair{u_r}{v_t}+\pair{u_t}{v_r}&=\beta_{rt}
 \qquad(r\ne t).
 \end{aligned}
\end{equation}
No division by two occurs. The coordinate recovery formula is
\begin{equation}\label{eq:char-two-recovery}
 s_rI_2=x_{r+2}+[E_{11},x_{r+2}],\qquad
 (u_r,v_r)=[E_{11},x_{r+2}].
\end{equation}
Together with the matrix units, it proves that $D$ is generated by the $s_r$ and all pairings of the $u_r,v_r$, and that these module vectors generate $M$. The proof is the same subalgebra argument as in Lemma~\ref{lem:universal-coordinates}.

\begin{proposition}\label{prop:char-two-decomposition}
If $\operatorname{rank}\beta\ge4$, the associative and octonion factors in Theorem~\ref{thm:char-two-main} exist. More explicitly, choose Clifford generators $x_1,x_2,x_3,x_4$ for two orthogonal nonsingular binary subforms, normalized by
\[
 \beta(x_1,x_2)=\beta(x_3,x_4)=1,
\]
where the notation for $\beta$ refers to the underlying generating vectors. Put
\begin{equation}\label{eq:char-two-intrinsic-e}
 P_r=x_r+[x_1x_2,x_r]\quad(r=3,4),\qquad e=[P_3,P_4].
\end{equation}
Then $e$ is the central idempotent with $eA(q)\simeq C(q)$.
\end{proposition}
\begin{proof}
First use the split coordinates above. Let $J$ be the pairing ideal in $D$, set
\[
 c_{rt}=[s_r,s_t],\qquad
 p=\pair{u_1}{v_2}+\pair{u_2}{v_1},\qquad c=c_{12}.
\]
Every $c_{rt}$ is central in $D$ by \eqref{eq:char-two-relations}. Since $[D,D]M=[D,D]J=0$, we have $cM=cJ=0$. The relation $\beta_{12}=1$ gives $c=1+p$, and therefore
\begin{equation}\label{eq:char-two-idempotent}
 c^2=c(1+p)=c.
\end{equation}
Lemma~\ref{lem:center} makes $e=cI_2$ a central idempotent of the whole alternative algebra.

In the $e$-component, $eM=eJ=0$, so $eA(q)=\Mat_2(eD)$ is associative. In the ordinary Clifford quotient the recovery formula gives $\pi(s_rI_2)=\pi(x_{r+2})$, and hence $\pi(e)=\beta_{12}=1$. Universality and generation give inverse homomorphisms $C(q)\to eA(q)\to C(q)$, just as in Proposition~\ref{prop:idempotent}.

Let $f=1-e$. In the coefficient algebra, $f=p\in J$. Thus
\[
 f[s_r,s_t]=0\quad\hbox{for all }r,t,
\]
and the generation statement following \eqref{eq:char-two-recovery} shows that $R=fD$ is commutative and finitely generated. Since $f\in fJ$, the pairing ideal is all of $R$. Lemma~\ref{lem:local-octonion} identifies $fA(q)$, if nonzero, as an octonion algebra with center $R$.

For descent, use the original normalized symplectic pairs $x_1,x_2$ and $x_3,x_4$, without requiring their binary forms to be hyperbolic over $F$. After a finite separable splitting extension, the element $h=x_1x_2$ in the quaternion subalgebra has trace one. The matrix component of each $x_r$ with $r\ge3$ is scalar, and
\[
 [h,sI_2+(u,v)]=(u,v)
\]
because $h+h^*=I_2$. Thus the $P_r$ in \eqref{eq:char-two-intrinsic-e} are exactly the scalar coordinate matrices, and their commutator is the idempotent just constructed. Formula~\eqref{eq:char-two-intrinsic-e} is defined over $F$, so its centrality and idempotence descend, as does $eA(q)\simeq C(q)$.

Centers commute with the finite field extension by comparison of coefficients in commutators and associators. Finite generation of the center and local freeness of rank eight descend as in Section~\ref{sec:decomposition}. The anti-involution fixing each Clifford generator restricts, over the splitting extension, to octonion conjugation: all the generators in \eqref{eq:char-two-generators}, and $i,j$, have trace zero. Thus the multiplicative nonsingular norm descends as well. Finally, $fA(q)=\ker\pi$ and $eA(q)$ is its two-sided annihilator, proving that the decomposition is canonical.
\end{proof}

\begin{corollary}\label{cor:char-two-cutoff}
If $\operatorname{rank}b_q\ge8$, then $A(q)=C(q)$. In particular, this holds for every nonsingular form of dimension at least eight.
\end{corollary}
\begin{proof}
It suffices to work after the splitting extension above. If the octonion factor were nonzero, take a residue field $L$ of its center. Its module $M_L$ is symplectic of dimension two. The images of the remaining generators define a linear map from the space of $q'_L$ to
\[
 W=L\oplus M_L\oplus M_L,\qquad
 Q(s,u,v)=s^2+\pair{u}{v},
\]
whose pullback quadratic form is $q'_L$, by \eqref{eq:char-two-relations} with commuting even coordinates. The polar form of $Q$ has rank four. Hence
\[
 \operatorname{rank}b_q-2
 =\operatorname{rank}b_{q'}\le4,
\]
a contradiction. This proves that the full octonion factor vanishes.
\end{proof}

Equivalently, the trace-zero space of an octonion algebra in characteristic two has dimension seven, with polar radical equal to the scalar line and polar rank six. This explains why polar rank, rather than dimension alone, controls the cutoff for possibly singular forms.

\subsection{The universal center and the four-dimensional quadric}
The normalized presentation also survives, with different quadratic equations. Retain $q=H\perp q'$ and $\beta_{12}=1$. Introduce $s_r$ and $p_{\alpha\beta}$ as in Section~\ref{sec:presentation}, now with $p_{\beta\alpha}=p_{\alpha\beta}$ and $p_{\alpha\alpha}=0$. Let $\mathcal R_2$ be the commutative ring defined by
\begin{equation}\label{eq:char-two-center-relations}
 \begin{aligned}
 p_{\alpha\beta}p_{\gamma\delta}
 +p_{\alpha\gamma}p_{\beta\delta}
 +p_{\alpha\delta}p_{\beta\gamma}&=0,\\
 p_{2r-1,2r}&=a_r+s_r^2,\\
 p_{2r-1,2t}+p_{2t-1,2r}&=\beta_{rt}\quad(r\ne t).
 \end{aligned}
\end{equation}
Let $P$ be the module defined by \eqref{eq:P-relations}, with pairing $\pair{z_\alpha}{z_\beta}=p_{\alpha\beta}$.

\begin{proposition}\label{prop:char-two-universal}
The pairings on $P$ are well defined and generate the unit ideal. Unless $\mathcal R_2=0$, the module $P$ is symplectic of rank two, and
\[
 R_q\simeq\mathcal R_2,\qquad
 \OO_q\simeq\mathcal A(\mathcal R_2,P).
\]
\end{proposition}
\begin{proof}
The Pl\"ucker relations make the pairing well defined. The equation
\[
 p_{14}+p_{23}=1
\]
shows that its values generate the unit ideal, so Lemma~\ref{lem:local-octonion} applies. In $\mathcal A(\mathcal R_2,P)$, the elements $E_{12},E_{21}$ and
\[
 s_rI_2+(z_{2r-1},z_{2r})
\]
satisfy the defining relations of $A(q)$ and send $e$ to zero. Conversely, the coordinates of the octonion factor satisfy \eqref{eq:char-two-center-relations}, giving a compatible ring and module map in the opposite direction. The two induced algebra homomorphisms are inverse by \eqref{eq:char-two-recovery} and generation by the matrix units and module vectors. This also identifies the centers and handles the zero-ring case.
\end{proof}

\begin{proposition}\label{prop:char-two-quadric}
Write $[a,b](x,y)=ax^2+xy+by^2$. If $q=H\perp[a,b]$, then
\begin{equation}\label{eq:char-two-quadric}
 R_q\simeq
 \frac{F[s,t,u,v,z]}
 {\bigl(uv+z^2+z+bs^2+st+at^2+ab\bigr)}.
\end{equation}
This is a smooth geometrically integral affine quadric of dimension four, and $\Cent(A(q))\simeq F\times R_q$.
\end{proposition}
\begin{proof}
There are four module vectors. Set
\[
 p_{12}=s^2+a,\quad p_{34}=t^2+b,\quad
 p_{13}=u,\quad p_{24}=v,\quad
 p_{14}=h,\quad p_{23}=1+h.
\]
The single independent Pl\"ucker relation gives
\[
 R_q\simeq
 F[s,t,u,v,h]/\bigl(uv+h^2+h+(s^2+a)(t^2+b)\bigr).
\]
Substitute $z=h+st$ to obtain \eqref{eq:char-two-quadric}. Its derivative with respect to $z$ is one, proving smoothness. Over an algebraic closure, its degree-two homogeneous part has polar rank four, whereas the degree-two part of a product of two affine linear polynomials has polar rank at most two. Thus the defining polynomial is irreducible. The center of the nonsingular four-dimensional ordinary Clifford algebra is $F$.
\end{proof}

For an arbitrary nonsingular four-dimensional $q$, the presentation \eqref{eq:char-two-quadric} holds after a finite separable extension making a chosen binary summand hyperbolic. The decomposition and its center are already defined over $F$ by Proposition~\ref{prop:char-two-decomposition}; formula~\eqref{eq:char-two-quadric} is asserted over $F$ itself only under the displayed hypothesis $q=H\perp[a,b]$.

\subsection{Geometry for nonsingular forms}
\begin{proposition}\label{prop:char-two-geometry}
For a nonsingular form $q$ of dimension $n=4$ or $6$, the ring $R_q$ is nonzero, smooth, and geometrically integral, of dimension $4$ or $7$, respectively.
\end{proposition}
\begin{proof}
Work over an algebraic closure $K$ and write $q=H\perp q'$. Put $m=n-2\in\{2,4\}$, let $U=K^2$ with its standard alternating form $\omega$, and set
\[
 T=U\oplus U,\qquad Q_T(u,v)=\omega(u,v),\qquad
 W=K\oplus T,
\]
with quadratic form $Q(s,y)=s^2+Q_T(y)$. The polar form $B$ of $Q_T$ is nonsingular on $T$. Let $X_m$ be the scheme of tuples $(s_r,y_r)_{r=1}^m$ with
\begin{equation}\label{eq:char-two-frame-equations}
 Q_T(y_r)+s_r^2=a_r,\qquad
 B(y_r,y_t)=\beta_{rt}\quad(r<t).
\end{equation}
By the same functorial argument as for \eqref{eq:sl2-torsor}, the map that records the $s_r$ and pairings of the module coordinates identifies
\[
 X_m\longrightarrow\Spec\mathcal R_2
\]
with the symplectic frame bundle of $P$. It is a Zariski locally trivial $\mathrm{SL}_2$-torsor: the opens where a $p_{\alpha\beta}$ is invertible cover the base, and the two corresponding vectors give a symplectic basis after rescaling.

We verify the geometry of $X_m$ without using the preceding orthogonal-group argument. Forgetting the $s_r$ maps $X_m$ to the symplectic frame variety
\[
 S_m=\{(y_r):B(y_r,y_t)=\beta_{rt}\},
\]
whose Gram matrix is nonsingular. After a choice of reference frame,
\[
 S_m\simeq\mathrm{Sp}_4/\mathrm{Sp}_{4-m}.
\]
It is nonempty and geometrically integral. The map $X_m\to S_m$ adjoins independent roots of the monic equations $s_r^2=a_r+Q_T(y_r)$. It is finite, faithfully flat, and radicial, hence a universal homeomorphism. Thus $X_m$ is nonempty and irreducible.

For smoothness, at any geometric point put
\[
 A_{rt}=B(y_r,\dot y_t).
\]
Varying the $\dot y_t$ in the span of the $y_r$ realizes every $m\times m$ matrix $A$, because $(\beta_{rt})$ is invertible. The differentials of the equations in \eqref{eq:char-two-frame-equations} are the independent linear forms
\[
 A_{rr}\quad(1\le r\le m),\qquad
 A_{rt}+A_{tr}\quad(r<t).
\]
The Jacobian criterion gives smoothness and
\[
 \dim X_m=5m-\frac{m(m+1)}2.
\]
Thus $X_m$ is integral. Smoothness descends along the $\mathrm{SL}_2$-torsor, and its surjectivity gives irreducibility of the base. Consequently $\mathcal R_2$ is nonzero, smooth, and integral over $K$, with
\[
 \dim\mathcal R_2=5m-\frac{m(m+1)}2-3
 =\begin{cases}4,&m=2,\\7,&m=4.\end{cases}
\]
These geometric assertions descend to $F$.
\end{proof}

\begin{corollary}\label{cor:char-two-semiprime}
If $q$ is nonsingular of dimension at least four in characteristic two, then $A(q)$ is semiprime and its center is reduced.
\end{corollary}
\begin{proof}
The ordinary Clifford algebra of a nonsingular even-dimensional quadratic form is central simple. In dimensions four and six, $R_q$ is a domain and $\OO_q$ is projective over it. Hence $\OO_q$ embeds into its scalar extension to $\Frac(R_q)$, which is a simple octonion algebra. A nilpotent ideal extends to zero there and therefore was zero. The central decomposition proves the claim. For dimension at least eight, Corollary~\ref{cor:char-two-cutoff} applies.
\end{proof}

Proposition~\ref{prop:char-two-decomposition}, Corollary~\ref{cor:char-two-cutoff}, Proposition~\ref{prop:char-two-geometry}, and Corollary~\ref{cor:char-two-semiprime} prove Theorem~\ref{thm:char-two-main}.

\subsection{Regular seven-dimensional forms and central nilpotents}
An odd-dimensional quadratic form in characteristic two is called regular here when its polar radical is one-dimensional and the quadratic form is nonzero on that line. Thus a regular seven-dimensional form can be written as $q_0\perp\langle a\rangle$ with $q_0$ nonsingular of dimension six and $a\ne0$.

\begin{proposition}\label{prop:char-two-seven}
Let $q=q_0\perp\langle a\rangle$ with $q_0$ nonsingular of dimension six and $a\ne0$. Put $L_a=F[t]/(t^2+a)$. Then
\begin{equation}\label{eq:char-two-seven}
 R_q\simeq R_{q_0}\otimes_F L_a,\qquad
 \OO_q\simeq\OO_{q_0}\otimes_F L_a.
\end{equation}
If $a=b^2$ in $F$, then $R_q$ contains a nonzero element $\varepsilon=t+b$ with $\varepsilon^2=0$, and $A(q)$ has a nonzero square-zero ideal. For every $a\ne0$, $R_q$ is geometrically nonreduced and has one geometric connected component, of dimension seven.
\end{proposition}
\begin{proof}
Choose the first two nonsingular binary summands inside $q_0$. Formula~\eqref{eq:char-two-intrinsic-e} shows that $A(q_0)\to A(q)$ preserves the central idempotents. We may check the remaining coordinate argument after a finite separable extension with $q_0=H\perp q'_0$.

In the octonion factor of $A(q)$, write the four generators of $q'_0$ as $s_rI_2+(u_r,v_r)$. The four vectors $(u_r,v_r)$ in the rank-four module $M\oplus M$ have nonsingular Gram matrix $(\beta_{rt})$. They therefore form a basis: the Gram matrix gives a splitting, and an isomorphism follows locally between modules of equal rank. The last generator, orthogonal to $q_0$, has the form
\[
 tI_2+(u,v).
\]
Its commutation relations with those four generators say that $(u,v)$ is orthogonal to that basis. Thus $(u,v)=0$. The last generator is consequently a scalar, and its square relation is $t^2=a$.

All remaining center and module relations are exactly those for $q_0$. Conversely, adjoining a central scalar $t$ with $t^2=a$ to $\OO_{q_0}$ supplies all the required Clifford relations. These constructions give inverse maps by the universal presentation of Proposition~\ref{prop:char-two-universal}. The scalar $t$ is intrinsically the image of the generator of $\langle a\rangle$, so the isomorphisms descend and give \eqref{eq:char-two-seven} over $F$.

If $a=b^2$, then $L_a=F[\varepsilon]/(\varepsilon^2)$. Since $R_{q_0}\ne0$, the element $\varepsilon$ is nonzero in $R_q$; the ideal $\varepsilon\OO_q$ is nonzero and has square zero. As $\OO_q$ is a direct factor of $A(q)$, the latter algebra is not semiprime. Over an algebraic closure, every $a$ is a square, and
\[
 R_q\otimes_F\overline F
 \simeq(R_{q_0}\otimes_F\overline F)[\varepsilon]/(\varepsilon^2).
\]
Proposition~\ref{prop:char-two-geometry} shows that its reduction is integral of dimension seven, proving the last assertion.
\end{proof}

For example, $q=H\perp H\perp H\perp\langle1\rangle$ has a nonzero central square-zero element in its octonion factor. This phenomenon belongs to characteristic two and does not affect the comparison in Section~\ref{sec:comparison}. If $a$ is not a square, $L_a$ is a purely inseparable field extension; the displayed center is then integral, although it is still geometrically nonreduced.

Finally, the rank cutoff still implies $A(q)=C(q)$ for every regular form of dimension at least eight: a nonsingular even-dimensional form has full polar rank, and a regular odd-dimensional form of dimension at least nine has polar rank at least eight. No cutoff depending only on dimension holds for arbitrary singular forms. Indeed, adjoining any number of zero quadratic coordinates to a nonsingular four-dimensional form gives a form whose alternative Clifford algebra still surjects onto the original nonassociative algebra, by sending the extra generators to zero.

\section{Musgrave's octonion Clifford algebra}\label{sec:octonion-clifford}
Musgrave's definition uses representations into octonion algebras over the ground field, rather than a list of polynomial identities. More precisely, \cite[Definition~3.3]{Musgrave} defines
\begin{equation}\label{eq:musgrave-definition}
 \Coct_F(q)=A(q)/K_F(q),\qquad
 K_F(q)=\bigcap_{\varphi\in\mathcal H_F(q)}\ker\varphi,
\end{equation}
where $\mathcal H_F(q)$ consists of all unital $F$-algebra homomorphisms from $A(q)$ to octonion algebras over $F$. If this set is empty, the intersection is $A(q)$ and the quotient is zero. Musgrave works in characteristic different from two; in characteristic two we use the same definition as an extension of that construction.

Throughout this section, assume either that $\operatorname{char}F\ne2$ and $q$ is nondegenerate of dimension at least four, or that $\operatorname{char}F=2$ and $\operatorname{rank}b_q\ge4$. Thus the canonical decomposition
\[
 A(q)=C(q)\times\OO_q
\]
is available. When $\OO_q\ne0$, write $R=R_q$ and $X=\Spec R$.

\begin{lemma}\label{lem:octonion-representations}
Let $E/F$ be a field extension and $O$ an octonion algebra over $E$. Every $F$-algebra homomorphism $\varphi:A(q)\to O$ has $\varphi(e)=0$, and its image generates $O$ as an $E$-algebra. If $\OO_q\ne0$, such a map is equivalently a pair consisting of an $F$-algebra homomorphism $\chi:R\to E$ and an isomorphism
\[
 \OO_q\otimes_{R,\chi}E\simeq O.
\]
In particular, the distinct kernels of the maps in $\mathcal H_F(q)$ are in bijection with $X(F)$.
\end{lemma}
\begin{proof}
In characteristic different from two, choose a nondegenerate four-dimensional subform. The images of its orthogonal generators are nonscalar and trace zero, as in Lemma~\ref{lem:seven}, and they generate $O$ by Lemma~\ref{lem:frame-generation}.

In characteristic two, choose a nonsingular four-dimensional subform and make its first binary summand hyperbolic after a finite separable extension of $E$. In the coordinates of $O$, the other two generators have the form $s_rI_2+(u_r,v_r)$ with $s_r$ scalar and
\[
 \pair{u_1}{v_2}+\pair{u_2}{v_1}=1.
\]
The four module vectors therefore span the symplectic module of rank two. Formula~\eqref{eq:char-two-recovery} and the matrix units show that the four Clifford generators generate the octonion algebra. Generation descends to $E$ by faithful scalar extension.

Thus the $E$-linear extension of $\varphi$ is surjective in either characteristic. It sends the central idempotent $e$ into the center $E$ of $O$. If its image were one, the representation would factor through the associative algebra $eA(q)\otimes_F E$, forcing $O$ to be associative. Hence $\varphi(e)=0$, and the map factors through $\OO_q$. Its restriction to $R$ takes values in $E$, since it commutes and associates with an $E$-spanning image. The induced map $\OO_q\otimes_{R,\chi}E\to O$ is a surjection between eight-dimensional $E$-algebras and is therefore an isomorphism. Conversely, every such pair gives a representation.

For $E=F$, the kernel is consequently
\begin{equation}\label{eq:octonion-point-kernel}
 C(q)\times\mathfrak m_\chi\OO_q,\qquad
 \mathfrak m_\chi=\ker\chi.
\end{equation}
Every $F$-point supplies its octonion fiber as a target. Distinct points give distinct kernels, because $\OO_q$ is faithfully locally free over $R$ and
$\mathfrak m_\chi\OO_q\cap R=\mathfrak m_\chi$.
\end{proof}

\begin{theorem}\label{thm:octonion-quotient}
If $\OO_q=0$, then $\Coct_F(q)=0$. Otherwise, define the rational-point vanishing ideal
\[
 I_F(R)=\bigcap_{\chi\in X(F)}\mathfrak m_\chi,
\]
with $I_F(R)=R$ when $X(F)$ is empty. Then
\begin{equation}\label{eq:octonion-quotient}
 K_F(q)=C(q)\times I_F(R)\OO_q,
 \qquad
 \Coct_F(q)\simeq\OO_q\otimes_R R/I_F(R).
\end{equation}
If nonzero, this quotient is an octonion algebra over the reduced ring $R/I_F(R)$, which is its center. It is semiprime.
\end{theorem}
\begin{proof}
If $\OO_q=0$, Lemma~\ref{lem:octonion-representations} excludes all octonion-valued representations. Otherwise, intersect \eqref{eq:octonion-point-kernel} over $X(F)$. For any finite projective $R$-module $P$ and any family of ideals $I_\lambda$,
\[
 \bigcap_\lambda I_\lambda P
 =\left(\bigcap_\lambda I_\lambda\right)P.
\]
Indeed, realize $P$ as a direct summand of $R^N$; then $I_\lambda P=P\cap I_\lambda R^N$, and projection onto $P$ gives the asserted equality. Apply this to $P=\OO_q$ to obtain \eqref{eq:octonion-quotient}.

The ring $R/I_F(R)$ is reduced because $I_F(R)$ is an intersection of maximal ideals. Base change preserves the octonion structure and identifies the center with this base ring. Finally, the projective octonion module embeds into the product of its fibers over the fraction fields of the minimal-prime quotients of that reduced noetherian ring. These fibers are simple octonion algebras, so a nilpotent ideal must be zero, as in Corollary~\ref{cor:semiprime}.
\end{proof}

Geometrically, $\Spec(R/I_F(R))$ is the reduced Zariski closure of $X(F)$ in $X$. Thus the octonion Clifford algebra is the restriction of the universal octonion factor to that closure. In particular, $\Coct_F(q)=\OO_q$ whenever $R$ is reduced and $X(F)$ is Zariski dense. Equality need not hold over an arbitrary ground field, even though the octonion factor itself is nonzero. The theorem proves Corollary~\ref{cor:octonion-main} apart from its final assertion, which follows below.

\begin{corollary}\label{cor:octonion-fields}
The following consequences hold under the standing hypotheses of this section.
\begin{enumerate}
\item If $F$ is algebraically closed and $\OO_q\ne0$, then
\[
 \Coct_F(q)\simeq\OO_q\otimes_R R_{\mathrm{red}},
 \qquad R_{\mathrm{red}}=R/\sqrt{0}.
\]
For $\operatorname{char}F\ne2$ and $4\le n=\dim q\le7$, this is $\OO_q$, whose center has respective dimension $4,6,7,7$. Hence $\Coct_F(q)$ is infinite-dimensional over $F$. For $n\ge8$, it is zero.
\item If $F$ is algebraically closed of characteristic two and $q$ is nonsingular of dimension four or six, then $\Coct_F(q)=\OO_q$ and is infinite-dimensional over $F$. If $\operatorname{rank}b_q\ge8$, the quotient is zero over any $F$ of characteristic two.
\item If $F$ is finite, then $X(F)$ is finite and
\[
 \Coct_F(q)\simeq\prod_{\chi\in X(F)}
 \bigl(\OO_q\otimes_{R,\chi}F\bigr),
 \qquad \dim_F\Coct_F(q)=8\,\#X(F).
\]
This formula is interpreted as the zero algebra when there are no rational points or $\OO_q=0$. \revision{For example, let $F=\mathbb F_p$ with $p$ odd and let $q=\langle a,b,c,d\rangle$ be nondegenerate. In \eqref{eq:quadric} the determinant of the five-variable form and the prescribed value both equal $abcd$. Their product is a nonzero square, so the finite-field point count gives $\#X(\mathbb F_p)=p^4+p^2$. Hence $\dim_{\mathbb F_p}\Coct_{\mathbb F_p}(q)=8(p^4+p^2)$.}
\end{enumerate}
\end{corollary}
\begin{proof}
For an algebraically closed field, the Nullstellensatz gives $I_F(R)=\sqrt{0}$, because $R$ is finitely generated. Theorems~\ref{thm:main} and~\ref{thm:char-two-main} give reduced centers of the stated positive dimensions in the indicated cases. \revision{Since every finite-dimensional commutative $F$-algebra is Artinian, $R$ is infinite-dimensional over $F$.} An octonion algebra is faithful of rank eight over its center, so \revision{$\OO_q$ is then infinite-dimensional over $F$ as well.} The vanishing assertions follow from the cutoff results.

Over a finite field there are finitely many rational points. Their maximal ideals are pairwise comaximal, and the Chinese remainder theorem gives
\[
 R/I_F(R)\simeq\prod_{\chi\in X(F)}F.
\]
Base change of $\OO_q$ proves the product and dimension formulas.
\end{proof}

The seven-dimensional characteristic-two exception has a particularly simple effect on this quotient.

\begin{corollary}\label{cor:octonion-char-two-seven}
Let $\operatorname{char}F=2$ and $q=q_0\perp\langle a\rangle$ with $q_0$ nonsingular of dimension six and $a\ne0$. Then
\[
 \Coct_F(q)=
 \begin{cases}
  0,&a\notin F^2,\\
  \Coct_F(q_0),&a\in F^2,
 \end{cases}
\]
where the second line denotes a canonical isomorphism induced by inclusion of $q_0$. In particular, over an algebraically closed field it is $\OO_{q_0}$, and the square-zero central thickening of $\OO_q$ disappears.
\end{corollary}
\begin{proof}
By Proposition~\ref{prop:char-two-seven}, $R_q=R_{q_0}\otimes_F F[t]/(t^2+a)$. If $a$ is not a square, this ring has no $F$-point. If $a=b^2$, every $F$-point sends $t$ to the unique root $b$, and
\[
 I_F(R_q)=I_F(R_{q_0})R_q+(t+b)R_q.
\]
The corresponding formula for the octonion factor in \eqref{eq:char-two-seven}, followed by Theorem~\ref{thm:octonion-quotient}, proves the result.
\end{proof}

For clarity, one can also define a different quotient by allowing octonion targets over all extension fields of $F$. Denote it by $\Coct_{\mathrm{ext}}(q)$. Then, whenever $\OO_q\ne0$,
\begin{equation}\label{eq:octonion-extension-quotient}
 \Coct_{\mathrm{ext}}(q)
 \simeq\OO_q\otimes_R R_{\mathrm{red}}.
\end{equation}
Indeed, Lemma~\ref{lem:octonion-representations} applies over every extension field. For each prime ideal $\mathfrak p\subset R$, the field $\Frac(R/\mathfrak p)$ and the corresponding octonion fiber give a representation whose kernel on $\OO_q$ is $\mathfrak p\OO_q$: the projective module over $R/\mathfrak p$ embeds into its fraction-field extension. Intersecting these kernels gives $\sqrt{0}\,\OO_q$. Every field-valued representation kills this ideal, proving \eqref{eq:octonion-extension-quotient}. This extension-field version equals the full octonion factor whenever $R$ is reduced; it should be distinguished from Musgrave's ground-field quotient \eqref{eq:musgrave-definition}.

\section{Consequences for earlier statements}\label{sec:comparison}
The decomposition and the description of the octonion quotient have consequences for two earlier accounts. We record them separately from the structural proofs, specifying the statements and hypotheses involved.

\subsection{Semiprimeness and the 2018 announcement}
\revision{Theorem~\ref{thm:main} gives semiprimeness for every nondegenerate $q$ of dimension at least four in characteristic different from two. In particular, a central element $g$ with $g^2=0$ must vanish, since it generates a square-zero ideal.} The report \cite[pp.~1272--1273]{Vishne}, describing joint work of Chapman, Rosenbaum, and Vishne, makes two types of assertion relevant here. Its Theorem~8 lists the center dimensions of the octonionic quotient and states that this quotient vanishes for $\dim q>7$. These assertions agree with the results above. \revision{For the nondegenerate form $q=\langle\alpha_1,\alpha_2,\alpha_3,\alpha_4\rangle$ and the elements $c',c_i''$ of its Lemma~3, the report gives the following statements.}
\begin{quote}
\revision{\textbf{Proposition 4.} \emph{The nonzero element $g=\sum\alpha_i^{-1}{c_i''}^{2}-{c'}^{2}$ is central, and satisfies $g^2=0$.}}

\revision{\textbf{Corollary 5.} \emph{The alternative Clifford algebra in dimension $\ge4$ is not semiprime.}}
\end{quote}

There is also a direct check on the expression for $c'$ in \cite[Lemma~3(1)]{Vishne}. The expression printed there is
\begin{equation}\label{eq:printed-c}
 \begin{split}
 c'={}&x_1(x_2(x_3x_4))+2[x_1x_2,x_3x_4]\\
 &-3\bigl((x_1x_2)(x_3x_4)
          +x_3((x_1x_2)x_4)-x_4((x_1x_2)x_3)\bigr).
 \end{split}
\end{equation}
In the ordinary Clifford quotient, let $\omega=x_1x_2x_3x_4$. The first term maps to $\omega$, the commutator maps to zero, and the three terms inside the parentheses map respectively to $\omega,\omega,-\omega$. Therefore
\[
 \pi(c')=\omega-3(\omega+\omega-(-\omega))=-8\omega.
\]
The volume element $\omega$ anticommutes with every $x_i$ and is invertible. Hence $-8\omega$ is not central when $\operatorname{char}F\ne2$. The image of a central element under a surjective homomorphism is central, so \eqref{eq:printed-c} cannot be central as printed.

This calculation uses only the expression printed in the announcement. \revision{A misprint could account for its failure to be central.} The comparison with its semiprimeness assertions follows independently from the decomposition proved here. The assertion of semiprimeness should also be distinguished from the absence of individual square-zero elements: split matrix and octonion algebras have such elements, but no nonzero nilpotent ideals.

\subsection{Finite dimensionality of the octonion Clifford algebra}
\revision{Chapman and Vishne \cite[Corollaries~4.4--4.5]{CV} showed that for a nondegenerate ternary form $q$ the center of $A(q)$ is $F[\delta]$. If $F$ is infinite of characteristic different from two, all but finitely many specializations $\delta\mapsto\lambda\in F$ give octonion fibers. Hence no nonzero polynomial in $\delta$ vanishes in the ground-field quotient $\Coct_F(q)$, already giving a negative answer to \cite[Question~3.4]{Musgrave} in dimension three. Corollary~\ref{cor:octonion-fields} extends this conclusion to $4\le n\le7$ over algebraically closed fields of characteristic different from two for Musgrave's ground-field quotient \eqref{eq:musgrave-definition}.}

\revision{For an octonion algebra over $F$ with norm form $\rho$, Musgrave writes in \cite[p.~590, Proposition~3.5]{Musgrave}, ``If $q$ [is] the pure part of the Pfister form $\rho$, then $Cl^{Oct}(q)$ is finite dimensional.''} For example, let $F=\mathbb C$, let $O$ be the split octonion algebra, and let $q$ be its seven-dimensional pure norm form. Over $F$, this form is isometric to its negative, so it satisfies both the pure-norm formulation of that proposition and the Clifford sign convention used here. Its octonion Clifford algebra has center of dimension seven and is infinite-dimensional over $F$.

The relevant step in that proof infers equality of kernels from isomorphism of the octonion targets. Isomorphic targets alone do not give equality of kernels. Lemma~\ref{lem:octonion-representations} instead identifies distinct kernels with distinct rational points of $X$. An isometry of an octonion norm fixing the identity need not preserve multiplication: over an algebraically closed field of characteristic different from two, the former group is $\mathrm O_7$, of dimension twenty-one, while the automorphism group is $G_2$, of dimension fourteen. The distinction is exactly the quotient $\mathrm O_7/G_2$ appearing in Proposition~\ref{prop:g2}.

\section*{Statements and Declarations}
\paragraph{Data availability.}
This article is theoretical. No research datasets were generated or analysed.

\paragraph{Use of generative artificial intelligence.}
The author consulted ChatGPT pro in the formulation of the proofs and in the preparation of this manuscript. The research questions and ideas came from the author, who takes full responsibility for the validity and coherence of this paper.


\end{document}